\documentclass[aps,pre,twocolumn,groupedaddress,amsmath,amssymb]{revtex4-1}
\usepackage{graphicx,subfigure,bm}

\begin{document}


\title{Solution of the Voter Model by Spectral Analysis}


\author{William Pickering}
\author{Chjan Lim}
\affiliation{Department of Mathematical Sciences, Rensselaer Polytechnic Institute, 110 8th Street, Troy, New York 12180, USA}

\date{\today}
\begin{abstract}
An exact spectral analysis of the Markov Propagator for the Voter
model is presented for the complete graph, and extended to the
complete bipartite graph and uncorrelated random networks. Using a
well-defined Martingale approximation in diffusion-dominated regions
of phase space, which is almost everywhere for the Voter model, this
method is applied to compute analytically several key quantities
such as exact expressions for the $m$ time step propagator of the
Voter model, all moments of consensus times, and the local times for
each macrostate. This spectral method is motivated by a related
method for solving the Ehrenfest Urn problem and by formulating the
Voter model on the complete graph as an Urn model. Comparisons of
the analytical results from the spectral method and numerical
results from Monte-Carlo simulations are presented to validate the
spectral method.

\end{abstract}


\maketitle

\section{Introduction}
\label{sec:intro} Quantitative insight into social opinion dynamics
can be found from a modeling perspective. The range of stochastic
models for opinion dynamics can be conveniently viewed as comprising
two main types: (I) Diffusion-dominated (II) Deterministic drift
dominated. One of the most well studied of the former is the Voter
model \cite{liggett,clifford}. This model has been solved
analytically for several quantities
\cite{vazquez,ben-naim,vilone,frachebourg,cox}, and is known to be a
perfect Martingale system on a complete graph. The second class
consists largely of stochastic models for social opinion dynamics
where the deterministic drift terms are significant for all
population sizes while the noise terms decreases with population
sizes. Important examples of the second class of models are the
Naming games \cite{xie,baronchelli2}.

We focus on formulating the Spectral method and using it to solve
the Voter model in this paper. Besides opinion dynamics, it also has
applications in many areas of physical science, such as particle
interactions and kinetics of catalytic reactions \cite{frachebourg}.
However, social opinion dynamics is a typical view of this particle
interaction model, where there are two opposing opinions in a
political discussion. Another interpretation of the voter model is that
species compete for territories, where a change in state corresponds
to and invasion \cite{clifford,liggett2}. The fundamental social assumption about these
social models is that individuals are strongly influenced by the beliefs of
their neighbors. This is particularly noticeable when there is a fraction of committed agents in the system. For some social models, a tipping point is observed in minority opinion over a critical value of these individuals \cite{xie,turalska}. Other extensions and
generalizations of the Voter model also have been studied in detail
\cite{sood, masuda, masuda2,caccioli,rogers,baronchelli}.

A specific and direct application of the voter model is to population genetics.
In the context of a binary coded genetic algorithm,
the voter model on a complete graph can be interpreted as the crossover
operator of bit strings over a single bit. This is the exchange of bits once two parent strings have been selected \cite{goldberg}.
The voter model is an instance of the genetic algorithm that examines the case where selection is uniform (constant fitness),
only one child string is generated per crossover, and there are no mutations. A biased voter model, which introduces a fitness value
for each state, also has been studied \cite{sood}. The genetic algorithm, with or without mutations, is a subject of interest from
an analytical perspective \cite{holland,goldberg,vose}. The work presented here will provide some additional analytical insight into
unbiased genetic drift.

Many mathematical models for opinion dynamics, including the Voter
model, can be viewed as discrete time Markov chains. It is known
that the eigenvalues and eigenvectors of the transition matrix of
the Markov chain have a vital role in the dynamics of the model. In
general, the second largest eigenvalue is related to an approximate
rate of convergence to stochastic equilibrium and/or absorbing
states \cite{boyd}. For the Voter model, this is when all nodes have
the same state. All dynamics of the model halt entirely once this
network state is attained. On a connected graph, the expected time
to attain this unanimous state can be bounded in terms of the
spectral gap of the transition matrix \cite{cooper}. The estimates
for the expected time to consensus that we provide are tighter than
the bounds given in Cooper \textit{et al.} \cite{cooper}, however we
do not propose that our estimates hold for all connected networks.
In addition, if one had access to the complete spectral
decomposition of the transition matrix, then any future probability
distribution can be computed in closed form. In this paper, we offer
a method for finding a closed form solution of the spectral problem
that yields such a decomposition for the voter model on the complete
graph, complete bipartite graph, and uncorrelated degree
heterogeneous networks.

The Voter model on the complete graph can be cast as an urn problem
related to but distinct from the Ehrenfest model \cite{ehrenfest}.
The classical Ehrenfest Urn problem consists of two urns with $N$
balls divided amongst them. In a single time step, one ball is
chosen at random and is transferred to the other urn. This process
is repeated ad infinitum. In the urn formulation of the Voter model
on the complete graph, two balls are selected, one after the other.
After selection, both of the balls are placed in the urn from which
the second ball came. A closed form diagonalization of the Markov
transition matrix for the Ehrenfest model was found by Mark Kac in
1947 using similar techniques we will utilize here \cite{kac}. We
will propose a generalization of those techniques which allows us to
solve the spectral problem for the Voter model. This particular
solution of the Voter model will give exact analytical expressions
for several quantities of interest. In particular, we can estimate
the moments of the consensus time for each graph. For the complete
graph, we will find the solution of the spectral problem exactly,
which yield further results. One such result is the expected time
spent on each macrostate prior to reaching consensus.

The outline of the paper is as follows. Section \ref{sec:model} will
describe more of the details of the model and the associated random
walk on the complete graph. In section \ref{sec:spectral}, we shall
derive a straightforward procedure to exactly solve the
corresponding eigenvalue problem of the Markov transition matrix.
With such a solution, one can compute the probability distribution
at any future time in closed form. In section
\ref{sec:applications}, we use this solution to find two quantities
for the complete graph model: the expected time before reaching
consensus and the expected time spent at each macrostate. These
expressions are functions of the initial probability distribution,
which is kept arbitrary. Sections \ref{sec:bipartite} and
\ref{sec:networks} will consider the complete bipartite graphs and
the uncorrelated heterogeneous random networks by extending the
compete graph solution. For these networks, we simplify the model by
assuming that the system is diffusion dominated, since this is the
primary avenue towards consensus \cite{sood}.

\section{The 2-State Voter Model}
\label{sec:model} We begin by assuming that the Voter model is
imposed upon a complete graph of $N$ nodes. More general graphs have
been studied as well \cite{yildiz,sood,vazquez,suchecki}. Sections
\ref{sec:bipartite} and \ref{sec:networks} below shall apply this
theory to more general networks. Regardless of network topology,
each node is assigned one of two states, $A$ or $B$. In a single
time step, a node is chosen randomly and will assume the state of
one of its neighbors, also chosen randomly \cite{clifford}. In this
procedure, it is possible that the network state may remain
unchanged for several time steps.

Let $n_A(m)$ and $n_B(m)$ represent the total number of agents
taking opinion $A$ and $B$ respectively at discrete time $m$. Since
the total number of nodes must be conserved and all nodes must take
one of these two opinions, $n_A(m)+n_B(m)=N$. Since $N$ is a
constant, this will allow us to simplify the model to a random walk
in a single variable, say $n_A$.

Now let us formalize the problem as a random walk in $n_A$. We write
it as
\begin{equation}
n_A(m+1)=n_A(m)+\Delta n_A(m).\label{randomwalk}
\end{equation}
For a given time step $m$, $n_A$ is considered to be a known constant, and the random behavior is exhibited in $\Delta n_A$.
Since only a single node is updated per time step, $\Delta n_A$ only takes values from $\{-1,0,1\}$. Let $p_j=\frac{j(N-j)}{N(N-1)}\label{pj}$.
Then, from the definition of the model, the probabilities of taking these values are
\begin{equation}Pr\{\Delta n_A(m)=1\; |\; n_A(m)=j\}=p_j,\end{equation}
\begin{equation}Pr\{\Delta n_A(m)=-1\; |\; n_A(m)=j\}=p_j,\end{equation}
\begin{equation}Pr\{\Delta n_A(m)=0\; |\; n_A(m)=j\}=1-2p_j.\end{equation}
Thus, the random walk is determined given an initial condition for $n_A(0)$.

\subsection{Markov Propagator By Generating Functions}
In this section, we will outline a general procedure that provides a recurrence relation for the probability distribution of
general random walks. The strategy is to construct a sequence of generating functions for the probability distribution of the
random walk. This process can also be applied to other models, incomplete graphs, or to multiple dimensions. An advantage of
this approach is that it is highly generalizable in these ways and only requires minor modification to do so. The procedure
can be applied to find the probability distribution of either the microstates or the macrostates of the model. Here,
the macrostate approach will be utilized since the model is imposed on a complete graph.

To begin, represent the probability distribution in generating function form. Let $a_j^{(m)}=Pr\{n_A(m)=j|n_A(0)=n\}$.
We introduce a sequence of generating functions $R^{(m)}(x)=\sum_j a_j^{(m)}x^j$ and seek to find a relationship between $R^{(m+1)}(x)$ and $R^{(m)}(x)$. From the random walk form of the model, the probability generating function for $\Delta n_A$ is $D_j(x)=p_jx+(1-2p_j)+p_jx^{-1}$. We will make use of the following properties of generating functions in the derivations to follow:

\textit{Product rule:} If $X$ and $Y$ are integer random variables with probability generating functions $F(x)$ and $G(x)$ respectively,
then the generating function of $X+Y$ is $F(x)G(x)$.

\textit{Sum rule:} If the probability space is partitioned into $N$ events, each with generating function $F_j(x)$, then the generating function for the entire space is $\sum_{j=1}^N F_j(x)$ \cite{newman,newman2,bender}.

For time step $m$, suppose that $n_A(m)=j$. Note that $a_j^{(m)}x^j$ is the corresponding generating function for this event. Now, utilize the product rule in equation \eqref{randomwalk} to deduce that the probability generating function for time $m+1$ is $a_j^{(m)}x^jD_j(x)$ in the event that $n_A(m)=j$. By the sum rule, we have that
\begin{equation}R^{(m+1)}(x)=\sum_{j=0}^N a_j^{(m)}x^jD_j(x).\label{propagator}\end{equation}
This is the generating function form of the Markov propagator of the random walk. This can be easily generalized to other models simply by
specifying the appropriate expression for $D_j(x)$ in equation \eqref{propagator}.
For the given network topology, the sum is simple enough to collect terms and obtain an explicit equation for the propagator:
\begin{equation}a_j^{(m+1)}=p_{j-1}a_{j-1}^{(m)}+(1-2p_j)a_j^{(m)}+p_{j+1}a_{j+1}^{(m)}\label{mequation}.\end{equation}

While the procedure to find a relationship between $R^{(m+1)}$ and $R^{(m)}$ can be highly generalizable,
we pose a new power series that will be utilized directly to solve the formulation in equation \eqref{mequation}.
Let $Q^{(m)}(x,y)=\sum_j a_j^{(m)}x^jy^{N-j}$. Some properties of generating functions of this type are listed below:

\begin{enumerate}
\item Multiply $Q^{(m)}$ by x/y to shift $a_j^{(m)}\rightarrow a_{j-1}^{(m)}$.
\item Multiply $Q^{(m)}$ by y/x to shift $a_j^{(m)}\rightarrow a_{j+1}^{(m)}$.
\item The generating function for $p_ja_j^{(m)}$ is $\frac{xy}{N(N-1)}Q^{(m)}_{xy}$.
\end{enumerate}
Using these properties, we rewrite equation \eqref{mequation} as

\begin{equation}
(x-y)^2Q^{(m)}_{xy}=N(N-1)\Delta_{+m} Q^{(m)}.
\end{equation}
Here, $\Delta_{+m}$ is the forward difference operator in the discrete variable $m$.
The goal of the subsequent section is to solve this equation explicitly to find the probability distribution for
arbitrary time steps. To do so, we turn to the spectral problem.

\section{The Spectral Problem}
\label{sec:spectral} It is clear from equation \eqref{mequation}
that the future probability distribution can be expressed as a
tri-diagonal transition matrix multiplied by the probability vector
$\mathbf{a}^{(m)}=[a_j^{(m)}]_{j=1}^N$. Solving the spectral problem
provides a basis for the initial distribution, which allows future
time steps to be computed in closed form. We shall find the solution
to this problem here.
\subsection{Matrix Eigenvalue Problem}

For eigenvalue $\lambda$ with eigenvector $\mathbf{v}=[c_j]_{j=0}^N$, we have that the spectral problem for the propagator can be written as
\begin{equation} p_{j-1}c_{j-1}+(1-2p_j)c_j+p_{j+1}c_{j+1}=\lambda c_j. \end{equation}
The solution to the problem begins by defining a multivariate generating function for $c_j$ as
\begin{equation}
G(x,y)=\sum_jc_jx^jy^{N-j}. \label{G}
\end{equation}
We seek a differential equation for the eigenvalue problem. This procedure is similar to Kac's solution of the Ehrenfest model \cite{kac}.
Here, however, we will obtain a PDE instead of an ODE. Using the formulation in section \ref{sec:model}
gives the following partial differential equation for $G$:
\begin{equation}
(x-y)^2G_{xy}=N(N-1)(\lambda-1) G.
\end{equation}
Note that this equation has no boundary or initial conditions. To find the appropriate solutions to this equation, we require that it
must take the form specified in equation \eqref{G}, not all coefficients $c_j=0$, and that $c_j=0$ for $j<0$ and $j>N$.

To solve this, use the change of variables $u=x-y$ and $H(u,y)=G(x,y)$ to transform the partial differential equation into
\begin{equation}
u^2(H_{uy}-H_{uu})=N(N-1)(\lambda-1) H.
\end{equation}
Since this is a linear change of variables, we still expect solutions to be of the form $H(u,y)=\sum_jb_ju^jy^{N-j}$.
Substituting this into the new partial differential equation for $H$ and collecting like terms will give the difference equation for $b_j$:
\begin{equation}
(j-1)(N-j+1)b_{j-1}=(j(j-1)+N(N-1)(\lambda-1)) b_j.
\end{equation}
Recall that we had required $c_j=0$ for $j<0$ and $j>N$. This is also true for $b_j$ since we applied a linear transformation.
Therefore, this difference equation would suggest that every $b_j=0$ unless it has singular behavior for some value of $j\in \{0...N\}$,
say when $j=k$. At this singular point, set the coefficient of $b_k$ to 0. This allows the value of $b_k$ to be non-zero, and thus not all $b_j=0$.
This argument determines the set of eigenvalues, which are found to be
\begin{equation}
\lambda_k=1-\frac{k(k-1)}{N(N-1)},\; k=0...N.
\end{equation}

Now we can find the eigenvectors. Since the equation for $b_j$ has singular behavior at $b_k$, the value at this point is arbitrary.
This is expected, since any multiple of an eigenvector remains an eigenvector. Without any loss, let $b_k=1$.
Now we can ascertain the explicit solution for $b_j$ when $k<j$:
\begin{equation}
b_j=\prod_{i=k+1}^{j} \frac{(i-1)(N-i+1)}{N(N-1)(\lambda_k-1)+i(i-1)}.\label{evector}
\end{equation}
When $j<k$, we have $b_j=0$. Thus all values of $b_j$ are determined for a given eigenvalue $\lambda_k$. Note that in \eqref{evector},
the requirement that $b_j=0$ is satisfied when $j>N$.

Now, we use $b_j$ to find the components of the eigenvector, $c_j$. To do this, express $H(u,y)$ in the original $x,y$ variables to obtain
\begin{align}
G(x,y)&=H(u,y)\\
&=\sum_i b_i(x-y)^iy^{N-i}\\
&=\sum_{i=0}^N \sum_{j=0}^i (-1)^{i-j}b_i{i\choose j}x^jy^{N-j}\\
&=\sum_{j=0}^N \sum_{i=j}^N (-1)^{i-j}b_i{i\choose j}x^jy^{N-j}.\label{evector2}
\end{align}
Therefore $c_j=\sum_{i=j}^N (-1)^{i-j}b_i {i\choose j}$, and thus the spectral problem is solved in closed form.

\subsection{Differential Eigenvalue Problem}
\begin{figure}
\includegraphics[scale=.5]{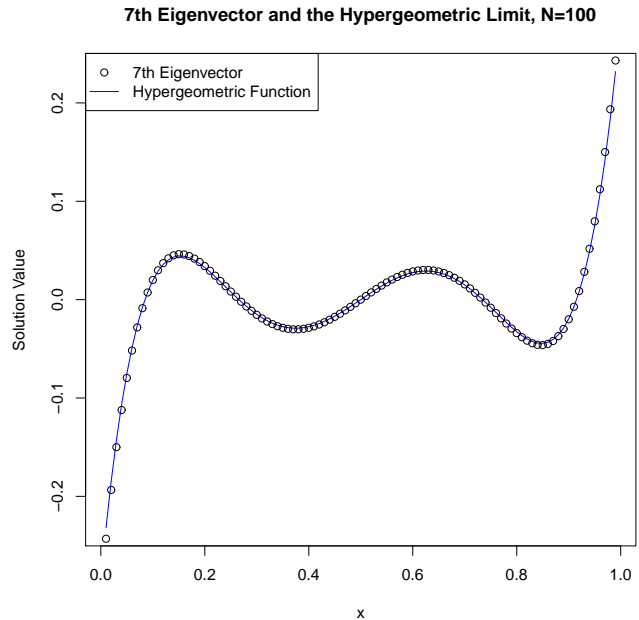}
\caption{7th Eigenvector of the discrete problem plotted with the exact solution for the limit as $N\rightarrow\infty$. The hypergeometric function in the figure is a fifth degree polynomial.\label{figure1}}
\end{figure}

In this section, we will examine the eigenvectors in more detail. In particular, we wish to consider the thermodynamic
limit of the model as $N\rightarrow\infty$ and study the behavior of the eigenvectors. To do this, note that the spectral problem can be posed as
\begin{equation}
\Delta^2_j(p_jc_j)=(\lambda_k-1) c_j,
\end{equation}
where $\Delta^2_j$ is the second centered difference operator over the discrete variable $j$. Let $x_j=j/N$, $\Delta x=1/N$ and $u(x_j)=c_j$. Then,
\begin{equation}
\frac{\Delta^2_j(p_ju(x_j))}{\Delta x^2}=N^2(\lambda_k-1) u(x_j).
\end{equation}
Thus, take $N\rightarrow\infty$ to obtain the differential equation for the eigenfunctions of the continuous time propagator:
\begin{equation}
\frac{d^2}{dx^2}[x(1-x)u(x)]=-k(k-1) u(x).
\end{equation}
Expanding all derivatives yields
\begin{equation}
x(1-x)\frac{d^2u}{dx^2}+(2-4x)\frac{du}{dx}+(k(k-1)-2)u=0,
\end{equation}
which is valid on $0<x<1$. This is a form of the hypergeometric differential equation, and as such, the basis for the solutions are
\begin{equation}
u_k(x)= {}_2F_1(k+1,2-k,2,x)\label{hypergeometric}
\end{equation}
for $k=0,1,2...$ This is a special case of the hypergeometric function in which the series expression terminates for each $k$ \cite{abramowitz}.
This implies that the \textit{kth} eigenvector of the voter model approaches a polynomial of degree $k-2$ as $N\rightarrow\infty$.
When $k=0$ and $k=1$, the solution is $u(x)=0$. This differential equation does not describe the behavior of the eigenvectors at
the boundary, which explains why the first two eigenfunctions are trivial. In figure \ref{figure1}, the seventh eigenvector
for the discrete model when $N=100$ has very close agreement with the continuous solution posed in \eqref{hypergeometric}.

\section{Applications of the Spectral Solution}
\label{sec:applications} The solution of the spectral problem yields
exact expressions for the probability distribution at any future
time step. With such a strong result, quantities that depend on
macro-state probabilities naturally follow from it. In this section,
we will apply the solution of the spectral problem to find two
quantities. The first is the expected time for the system to reach
consensus and the second is the expected time spent at each
macro-state before consensus, which is known as the local times.
Both quantities are considered functions of the initial distribution
of macro-states.

\subsection{Moments of Consensus Time}
A topic of particular importance with social models such as this is
the expected time until all agents in the network have the same
state. Once such a state is achieved, the dynamics halt entirely. In
this section, we will provide exact expressions for not only the
consensus time, but all moments thereof. This calculation depends on
the initial distribution of the model; however an estimate that is
independent of initial data shall also be determined. In deriving
the asymptotic bound that is independent of initial distribution,
the methods shall be left general enough to apply to any of the
networks we consider below. It is important to note that we use the
term \textit{time} to mean \textit{iterations} of the voter model
since the spectrum is utilized at the discrete level.

We start by finding the exact solution to the moments of consensus time on the complete graph.
We know that if $q_m$ is the probability of reaching consensus at step $m$,
then the $pth$ moment of consensus time is $\sum_{m=1}^\infty q_mm^p$.
We can determine $q_m$ from the boundary conditions of the transition equation.
The probability of reaching consensus at time step $m+1$ is the probability that the model
is in consensus at time $m+1$ minus the probability it was already in consensus at time $m$.
Symbolically, this is $q_{m+1}=a_0^{(m+1)}+a_N^{(m+1)}-a_0^{(m)}-a_N^{(m)}$. Now from the boundary conditions, we have that
\begin{equation}
a_{0,N}^{(m+1)}-a_{0,N}^{(m)}=\frac{1}{N}a_{1,N-1}^{(m)} \label{boundary}
\end{equation}
where we define $a_{0,N}^{(m)}$ to mean that this equation applies to either boundary $a_{0}^{(m)}$ or $a_{N}^{(m)}$ respectively.
Thus $q_m=(a_{1}^{(m-1)}+a_{N-1}^{(m-1)})/N$. Let $\mathbf{v_k}$ be the $kth$ eigenvector determined by combining equations \eqref{evector}
and \eqref{evector2} from section \ref{sec:spectral}. If we let $\mathbf{a^{(m)}}$ takes components $a_j^{(m)}$, then we can express
any such distribution as $\mathbf{a^{(m)}}=\sum_{k=0}^N d_k\lambda_k^m\mathbf{v_k}$. Here, $d_k$ is the initial probability distribution expressed
in the eigenbasis determined by solving the equation $\mathbf{Vd}=\mathbf{a^{(0)}}$, where $\mathbf{V}$ is the matrix whose columns are $\mathbf{v}_k$.
We only need components $2$ and $N-1$ from this expression, so we let $s_k=d_k\{[\mathbf{v_k}]_2+[\mathbf{v_k}]_{N-1}\}$ and write
\begin{equation}
a_{1}^{(m)}+a_{N-1}^{(m)}=\sum_{k=2}^N s_k\lambda_k^m.\label{bottleneck}
\end{equation}

We exclude the $k=0$ and $k=1$ terms since $a_1^{(m)}$ and $a_{N-1}^{(m)}$ are independent of those eigenvectors.
We now can write the $pth$ moment of the time to consensus as

\begin{align}
E[T^p|\mathbf{a^{(0)}}=\mathbf{Vd}]&=\sum_{m=1}^\infty \frac{1}{N}\sum_{k=2}^N s_k\lambda_k^{m-1}m^p\\
&\sim\sum_{k=2}^N \frac{1}{N}s_k \frac{p!}{(1-\lambda_k)^{p+1}}\\
&=p![N(N-1)]^{p+1}\sum_{k=2}^N  \frac{s_k}{[k(k-1)]^{p+1}}.
\end{align}

Now let us estimate this quantity to ascertain asymptotic information of the moments of consensus times.
The process outlined here will apply to general networks, not just the complete graph.
Let $S^{(m)}$ be the probability that the system is not in consensus at time $m$.
This is identical to the sum over probability of all non-consensus states.
As these interior probabilities are independent of the eigenvectors corresponding to equilibrium,
they are dominated by the second largest eigenvalue, $\lambda_2$. Since there is no drift, we have that $S^{(m)}= O(\lambda_2^m)$.

To find the moments, we need the probability of entering consensus at time $m$.
In this notation, we have that $q_m=S^{(m-1)}-S^{(m)}=O(\lambda_2^m(\frac{1-\lambda_2}{\lambda_2}))$.
Therefore, we have that the moments of consensus time can be bounded as follows:

\begin{align}
E[T^p]&=\sum_{m=0}^\infty q_mm^p\\
&=\frac{1-\lambda_2}{\lambda_2}\sum_{m=0}^\infty O(\lambda_2^m)m^p\\
&\sim \frac{1-\lambda_2}{\lambda_2}O\left(\frac{p!\lambda_2^{p+1}}{(1-\lambda_2)^{p+1}}\right)\\
&=O\left(\frac{p!\lambda_2^{p}}{(1-\lambda_2)^p}\right)\label{moments}
\end{align}
This is a uniform estimate for general initial probability
distributions for the $pth$ moment of consensus times of the model.
For the expected time to consensus $(p=1)$ on the complete graph the
asymptotic behavior of the expected time to consensus is quadratic
with $N$, which is consistent with previous results
\cite{sood,yildiz,zhang}. Using continuous time methods, the
expected number of iterations to reach consensus given a opinion
density $\rho$ is
$E[T]=N^2\left[(1-\rho)\ln\frac{1}{1-\rho}+\rho\ln\frac{1}{\rho}\right]$
\cite {sood}. The exact solution and the bound we found here are
improvements upon this result since it provides an expression for
all moments of consensus time, the solution valid even for small
values of $N$ or crowds, and the theory can be applied to general
initial probability distributions. Furthermore, the derivation can
be applied to any of the networks we shall consider.

\subsection{Local Times}
We define \textit{local times} as the expected frequency of
visitations of each macro-state of the random walk prior to
consensus. Unlike the consensus time, local times will be organized
as a vector whose components correspond to each state. The consensus
states are not included because are understood to have an infinite
local time. For non-consensus states, we expect local times to be
finite since the consensus time is also finite. Local times are more
detailed  quantities than the consensus time since the sum of all
local times is identical to the consensus time. We shall assume
throughout that the model is on the complete graph.

\subsubsection{Discrete Time Solution}

Since the solution of the spectral problem on the complete graph is known, computing local times becomes straightforward.
Let $M_j(m)$ be the total number of visitations of macrostate $j$ by time $m$ and let $\Delta M_j(m)=M_j(m)-M_j(m-1)$.
Note that $M_j(m)$ depends on the outcome of the random walk for $n_A$. Since $a_j^{(m)}$ is defined as the probability that $n_A(m)=j$, we have that $\Delta M_j(m)$ takes value 1 with probability $a_j^{(m)}$. Otherwise, it takes value 0. Thus, for $j=1...N-1$, we can write the local time as

\begin{equation}
E[M_j]=\sum_{m=0}^{\infty} E[\Delta M_j(m)]=\sum_{m=0}^{\infty} a_j^{(m)}.
\end{equation}

\begin{figure}
\includegraphics[scale=.4]{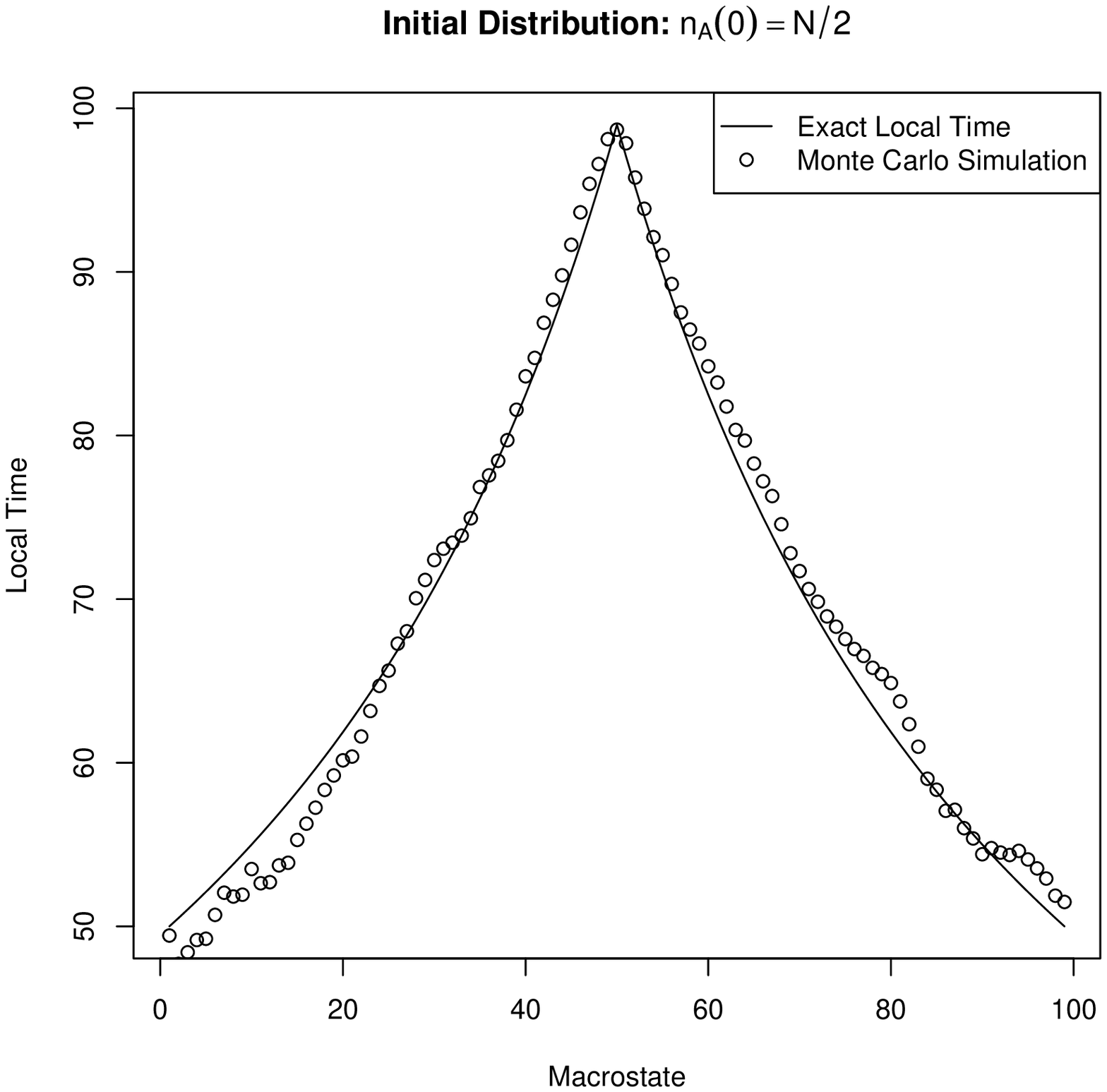}
\includegraphics[scale=.4]{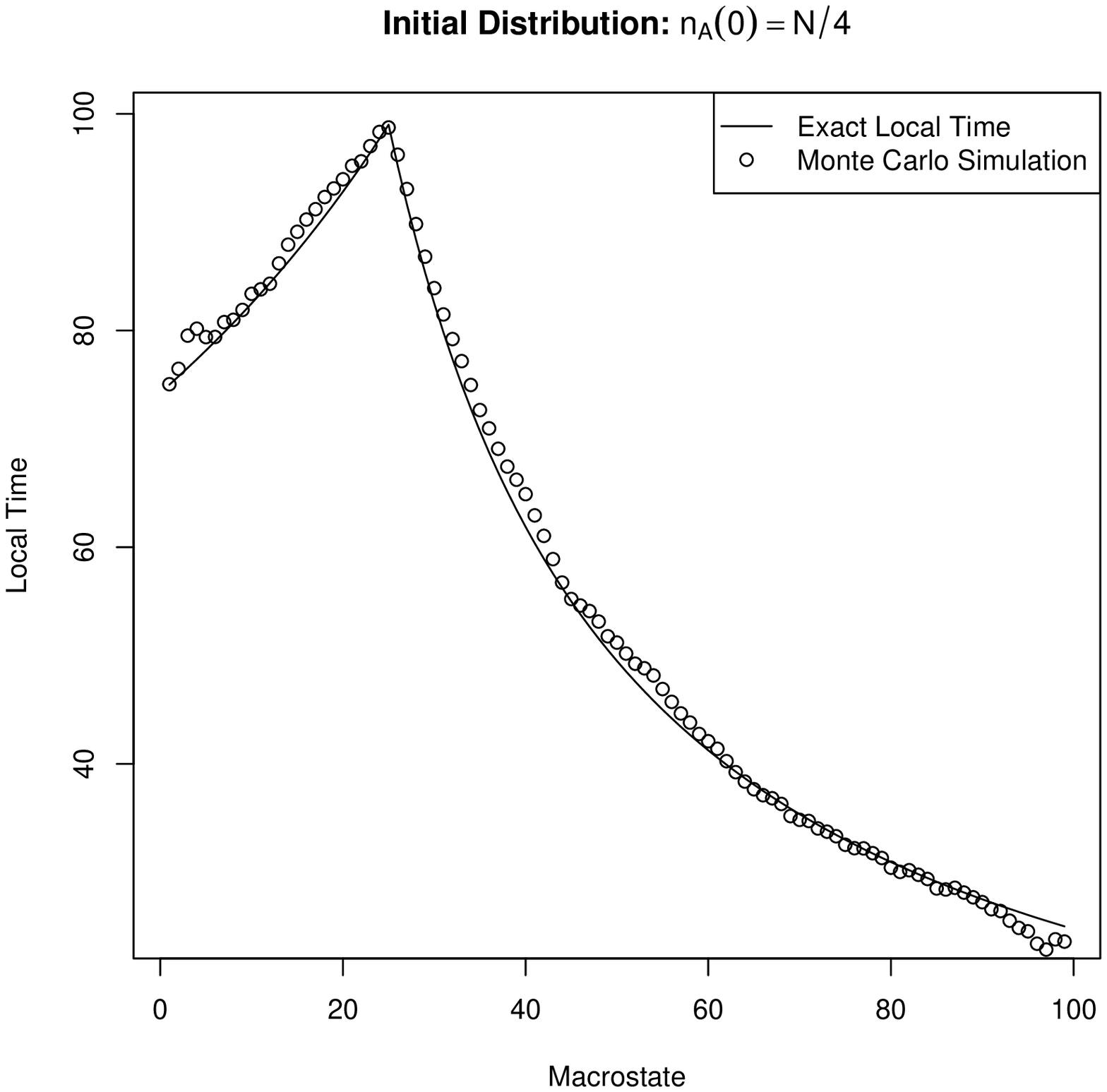}
\includegraphics[scale=.4]{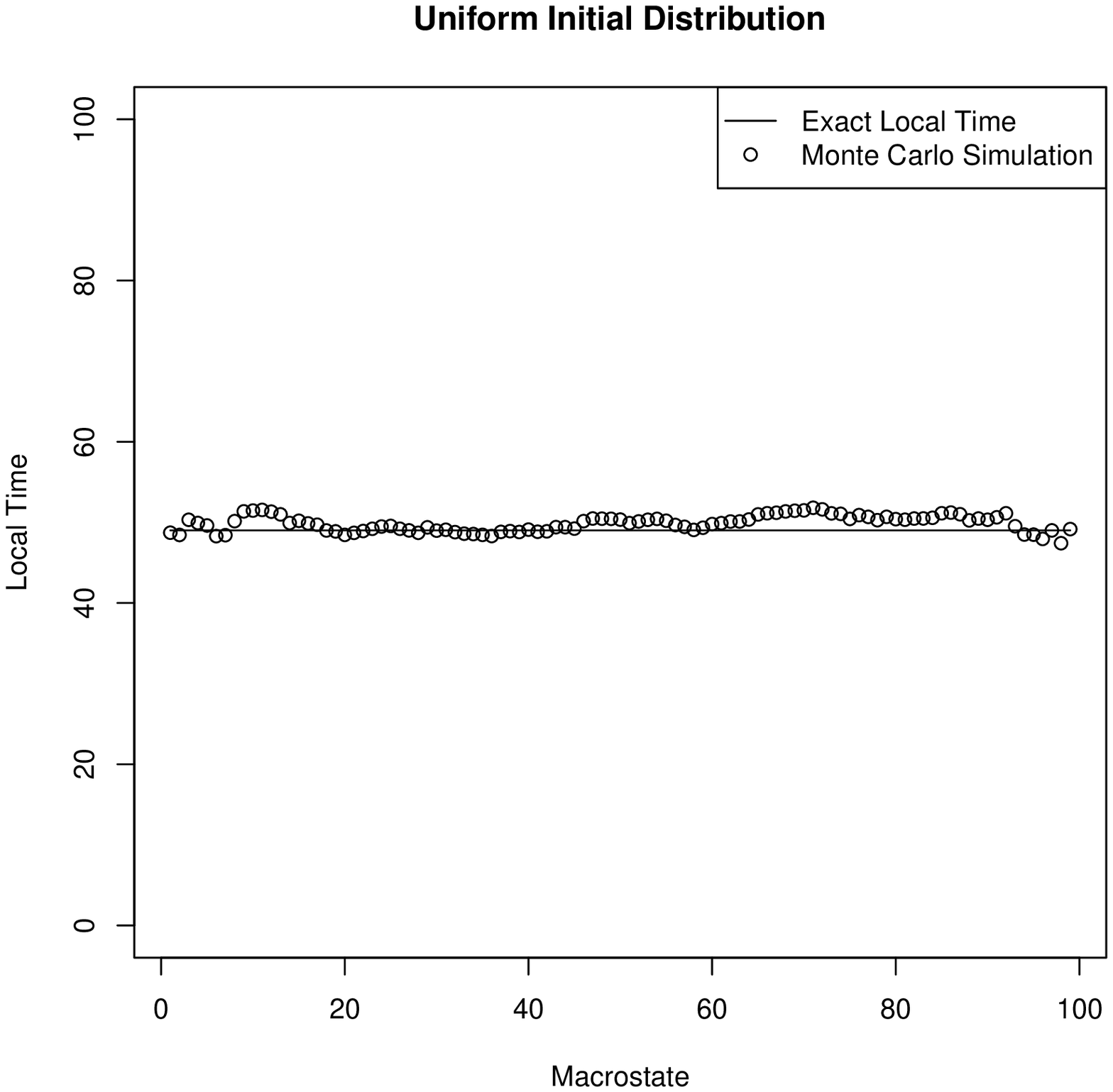}
\caption{Three cases of the initial distribution are examined. The exact expression for the local time in equation \eqref{macrostatetime} is compared with averaged results from a Monte Carlo simulation over 3000 runs of the voter model when $N=100$.\label{mcmc}}
\end{figure}

We can use the solution to the spectral problem to compute this infinite series exactly.
Let $\mathbf{M}$ be a vector whose components are $E[M_j]$. Also, letting $d_k$ be the initial distribution
of the macrostates expressed in the eigenbasis, we have that the time $m$ probability distribution
is $\mathbf{a^{(m)}}=\sum_{k=0}^{N} d_k\lambda_k^m\mathbf{v}_k$. The $k=0$ and $k=1$ terms in the sum
are the contributions of the consensus states to the probability distribution. Since consensus is a frozen state,
the probability distribution for the relevant macrostates are independent of the first two terms in this sum. The
local times can be found once these terms are discarded:

\begin{align}
\mathbf{M}&=\sum_{m=0}^{\infty}\sum_{k=2}^{N} d_k\mathbf{v}_k\lambda_k^m\\
&=N(N-1)\sum_{k=2}^{N} \frac{d_k}{k(k-1)}\mathbf{v}_k\label{macrostatetime}
\end{align}

The first and last components of this vector are irrelevant since it is understood that $M_0=M_N=\infty$.
The remaining components of $\mathbf{M}$ are the exact values for the local times.

We perform Monte Carlo simulation to reinforce this result. Three cases of the initial distribution are
considered: $a_j^{(0)}=\delta_{j,N/2}$, $a_j^{(0)}=\delta_{j,N/4}$, and $a_j^{(0)}=\frac{1}{N+1}\; \forall j=0...N$.
Figure \ref{mcmc} shows that there is good agreement between the exact solutions and the results from the simulations.
A particularly interesting case is when the initial distribution is uniform. This distribution happens to be the
eigenvector corresponding to the second largest eigenvalue of the transition matrix. So, in equation \eqref{macrostatetime},
we know that $d_2=\frac{1}{N+1}$ and $d_k=0$ otherwise. Therefore, the time spent at each macrostate is also uniform with
value $\frac{N(N-1)}{2(N+1)}\sim \frac{N}{2}$.

\subsubsection{Continuous Time Solution}

Here, we will find the expected time spent per macro-state as
$N\rightarrow\infty$ using continuous time methods. While the
discrete time solution is valid for every $N$, the continuous time
solution will provide some additional insight into its behavior for
large $N$. Let $\rho_j=j/N$, $t_m=m/N$, and
$u(\rho_j,t_m)=Na_j^{(m)}$. Using this, the Fokker-Plank equation
for the model is
\begin{equation}
\frac{\partial u}{\partial t}=\frac{1}{N}\frac{\partial^2}{\partial \rho^2}[\rho(1-\rho)u(\rho,t)].\label{fokkerplank}
\end{equation}
Furthermore, let $M(\rho)$ be the expected time spent at density $\rho$ prior to consensus. As $N\rightarrow \infty$, we have that
\begin{equation}
M(\rho)=\sum_{m=0}^\infty a_j^{(m)}\rightarrow\int_0^\infty u(\rho,t)dt.
\end{equation}
Now, integrate equation \eqref{fokkerplank} from $t=0$ to $t=\infty$ to obtain

\begin{equation}
u(\rho,\infty)-u(\rho,0)=\frac{1}{N}\frac{\partial^2}{\partial \rho^2}[\rho(1-\rho)M(\rho)].
\end{equation}
Since the system reaches consensus in finite time, we have that $u(\rho,\infty)=0$.
Let $f(\rho)=u(\rho,0)$ and $T(\rho)=\rho(1-\rho)M(\rho)/N$. Then, the problem becomes

\begin{equation}
-f(\rho)=\frac{\partial^2T}{\partial \rho^2}
\end{equation}
with boundary conditions $T(0)=T(1)=0$. The solution of this problem is found by determining the Green's function of
the differential operator \cite{herron}. For this problem, the Green's function is given by
\begin{equation}
\tilde{g}(\rho,\xi)=\left\{
\begin{array}{lr}
\rho(1-\xi), & \rho<\xi\\
\xi(1-\rho), & \rho>\xi\\
\end{array}
\right. .
\end{equation}
With this solution, we have that
\begin{equation}
T(\rho)=\int_0^1 f(\xi)\tilde{g}(\rho,\xi)d\xi.
\end{equation}
Therefore, the expression for local time for large $N$ is

\begin{equation}
M(\rho)\sim N\int_0^1 f(\xi)g(\rho,\xi)d\xi.
\end{equation}
where
\begin{equation}
g(\rho,\xi)=\left\{
\begin{array}{lr}
\frac{1-\xi}{1-\rho}, & \rho<\xi\\
\frac{\xi}{\rho}, & \rho>\xi\\
\end{array}
\right. .
\end{equation}

The Green's function $g(\rho,\xi)$ is the solution to local times when the initial density is specified as a
given value $\xi$. Furthermore, when the initial probability distribution is uniform, the solution is also uniform
with value $M(\rho)=N/2$ as we have observed with the discrete formulation.

\section{The Complete Bipartite Graph}
\label{sec:bipartite}
The methods we have developed for the complete graph can be extended to more complex networks. In particular,
we will consider the complete bipartite graph. In this case, nodes in the network are divided into two groups.
Every node in a group is connected to every node in the other group. The complete bipartite graph can also be
defined as the complement of two complete graphs. Let $N_1$ be the number of nodes in the first group and $N_2$
be the total number of nodes in the second group. Also, let $n_A^{(1)}(m)$ and $n_A^{(2)}(m)$ be the number of nodes with
opinion $A$ in groups 1 and 2 respectively.

Letting $a_{ij}^{(m)}=Pr\{n_A^{(1)}(m)=i,\; n_A^{(2)}(m)=j\}$ and $Q^{(m)}(x,y,u,v)=\sum_{i,j}a_{ij}^{(m)}x^i,y^{N_1-i}u^jv^{N_2-j}$, we can use the procedure laid out in section \ref{sec:model} to find the single step propagator for the model:

\begin{multline}
\left[\frac{u(x-y)}{NN_2}-\frac{y(u-v)}{NN_1}\right]Q_{yu}^{(m)}\\
+\left[\frac{-v(x-y)}{NN_2}+\frac{x(u-v)}{NN_1}\right]Q_{xv}^{(m)}=\Delta_{+m}Q^{(m)}.\label{bpropagator}
\end{multline}

Unlike for the complete graph, the closed form expression for the solution to the spectral problem for this equation will not be found with these methods.
However, we can apply some assumptions to the system that can reduce propagator to the region of the $(i,j)$ grid that is diffusion dominant.
Then, we wish to find the approximate size of the spectral gap, since this governs the expected time to consensus. Since we need to restrict
the region to estimate the spectrum, it will not allow us to find the eigenvectors exactly. As such, the more detailed solutions that depend on
the behavior of the eigenvectors, such as local times, will not be found.

If we take the single step propagator to continuous time, it is
known that the system approaches equilibrium when
$n_A^{(1)}/N_1-n_A^{(2)}/N_2\sim 0$ \cite{sood}. Furthermore, the
time to reach this equilibrium state is negligible compared to the
time to reach consensus. Along this line, diffusion governs the
motion of the macro-state of the system instead of drift. The study
of the behavior of the probability distribution is most valuable
when drift can be neglected. As such, we will make this assumption
to obtain

\begin{equation}
yuQ^{(m)}_{yu}\approx xvQ^{(m)}_{xv}.
\end{equation}
With this simplification, there are three cases that need to be accounted for:

\begin{enumerate}
\item $yu\approx xv \longrightarrow Q^{(m)}_{yu}\approx Q^{(m)}_{xv}$
\item $yu\ll xv \longrightarrow Q^{(m)}_{xv}\ll Q^{(m)}_{yu}$
\item $xv\ll yu \longrightarrow Q^{(m)}_{yu}\ll Q^{(m)}_{xv}$
\end{enumerate}

Consider the first case. Since $Q^{(m)}_{yu}\approx Q^{(m)}_{xv}$, we can combine the two terms in equation \eqref{bpropagator} together to obtain

\begin{equation}
\frac{(u-v)(x-y)}{N_1N_2}Q^{(m)}_{yu}\approx \Delta_{+m}Q^{(m)}.
\end{equation}
Letting $G(x,y,u,v)=\sum_{ij} c_{ij}x^iy^{N_1-i}u^jv^{N_2-j}$, the spectral problem is given by
\begin{equation}
(u-v)(x-y)G_{yu}\approx N_1N_2(\lambda-1)G.\label{simplified}
\end{equation}
Let $s=u-v$, $r=x-y$, and $H(r,y,s,v)=G(x,y,u,v)$. As with the complete graph, this is a linear transformation,
so we expect solutions of the same form: $H(r,y,s,v)=\sum_{ij} b_{ij}r^iy^{N_1-i}s^jv^{N_2-j}$. The differential equation becomes

\begin{equation}
rs(H_{sy}-H_{rs})\approx N_1N_2(\lambda-1)H.
\end{equation}
The corresponding difference equation for the coefficients, $b_{ij}$, is

\begin{equation}
j(N_1-i+1)b_{i-1,j}-ijb_{ij}\approx N_1N_2(\lambda-1)b_{ij}.\label{bipdiffeq}
\end{equation}
With similar arguments for the complete graph, we require a singularity in this difference equation
so that the solution is not trivial. Also, we assume that both sides of \eqref{bipdiffeq} agree to within an $O(1)$
factor since this is an approximation of the propagator. So, we write that $O(ij)=N_1N_2(\lambda-1)$. Taking $i,j$
to be small yields the approximate size of the spectral gap:

\begin{equation}
\lambda_{ij}= 1-O\left(\frac{1}{N_1N_2}\right).\label{eigen1}
\end{equation}

Note that in procedure, the same spectrum would be recovered if one were to use $G_{xv}$ instead of $G_{yu}$ in equation \eqref{simplified}.

Now let us consider cases 2 and 3 together. Without loss of generality, results from case 3 can be recovered by case 2 through
interchanging $u\leftrightarrow x$, $y\leftrightarrow v$, and $N_1\leftrightarrow N_2$. Physically, this corresponds to
interchanging the labels on the two groups of nodes. Taking case 2 as the archetype, drop small terms to simplify the propagator to

\begin{equation}
\frac{ux}{NN_2}Q^{(m)}_{yu}+\frac{yv}{NN_1}Q^{(m)}_{yu}-\frac{2yu}{N_1N_2}Q^{(m)}_{yu}\approx\Delta_{+m}Q^{(m)}.
\end{equation}

There is no need to change variables to solve the corresponding spectral problem since the corresponding finite difference
equation explicitly determines $c_{ij}$. Using similar arguments as in case 1, the spectral gap agrees for this case:
\begin{equation}
\lambda_{ij}=1-O\left(\frac{1}{N_1N_2}\right).\label{eigen2}
\end{equation}

Since this result is consistent for all regions of the $(x,y,u,v)$ plane that we consider, we have the asymptotic form of the spectral gap.

We can find the moments of consensus times by using equation \eqref{moments} from section \ref{sec:applications}. Taking $p=1$,
this suggests that the expected time to consensus is $O(N_1N_2)$. This is consistent with continuous time analysis, which shows
that the expected number iterations to reach consensus is $E[T]=4N_1N_2\left[(1-\omega)\ln\frac{1}{1-\omega}+\omega\ln\frac{1}{\omega}\right]$,
where $\omega$ is the degree weighted mean of microstates \cite{sood}. As before, this bound is valid for all initial probability distributions.
However, without detailed information about the eigenvectors, we cannot extract more detailed information about the propagator than the bound on
moments of consensus time.

\section{Uncorrelated Heterogeneous Networks}
\label{sec:networks} Here we will look into the spectrum when the
Voter model is imposed on a networks with arbitrary, fixed degree
sequences. The expected number of iterations to reach consensus on
these networks is known to be $O(N^2\mu_1^2/\mu_2)$, where $\mu_1$
and $\mu_2$ are the first and second moments of the degree sequence
respectively \cite{sood,vazquez}. We shall verify this result by
investigating the spectrum of the transition matrix, which we have
shown to give control over the expected time to consensus.
Furthermore, we can utilize this solution to find the moments of
consensus time.

To begin, we establish the following notation. Let $\mathbf{n}\in\mathbb{Z}_2^N$ be
a vector consisting of the microstates of each node in the network. Components of $\mathbf{n}$
take value $1$ when the node has state $A$ and takes value $0$ when it has state $B$. Let $\mathbf{A}$
be the adjacency matrix of the network, and $\mathbf{k}\in \mathbb{R}^N$ contain the degrees of each node.
Also, let component $k$ of $\mathbf{n}_A(m)$ be the number of nodes with degree $k$ that have opinion $A$.
Let $N_k$ be the total number of nodes with degree $k$. We also take $\{\mathbf{e}_i\}_{i=1}^N$ to be the
standard basis vectors for $\mathbb{R}^N$.

Let $a_{\bm{\alpha}}^{(m)}=Pr\{\mathbf{n}_A(m)=\bm{\alpha}\}$. Then, the single step propagator is found to be

\begin{multline}
\Delta_{+m}a_{\bm{\alpha}}^{(m)}=\sum_k [p_{k\bm{\alpha}-\bm{e}_k}a_{\bm{\alpha}-\bm{e}_k}^{(m)}-(p_{k\bm{\alpha}}+q_{k\bm{\alpha}})a_{\bm{\alpha}}^{(m)}\\
+q_{k\bm{\alpha}+\bm{e}_k}a_{\bm{\alpha}+\bm{e}_k}^{(m)}].
\end{multline}
Here, the transition probabilities are
\begin{align}
p_{k\bm{\alpha}}&=\frac{1}{Nk}\sum_{i\in K} \mathbf{Ae}_i\cdot\mathbf{n}(1-n_i)\\
q_{k\bm{\alpha}}&=\frac{1}{Nk}\sum_{i\in K} \mathbf{Ae}_i\cdot(\mathbf{1}-\mathbf{n})n_i,\\
\end{align}
where $K$ is the set of all nodes with degree $k$. As with the complete bipartite graph, let us assume that
the system is diffusion driven. The system will quickly approach this state in $O(N)$ steps and move towards
consensus on a much slower time scale \cite{sood}. Furthermore, we replace $\mathbf{A}$ with the mean over all
adjacency matrices to give an expected time to consensus over random networks. This suggests that we take
\begin{equation}
A_{ij}\rightarrow \frac{k_ik_j}{N\mu_1}.
\end{equation}
With these assumptions, we require that for each $k$,
\begin{align}
&\frac{\alpha_k}{N_k}\approx\frac{\mathbf{k\cdot n}}{N\mu_1}\label{alpharestriction}\\
&q_{k\bm{\alpha}}\sim p_{k\bm{\alpha}}\sim \frac{\alpha_k(N_k-\alpha_k)}{NN_k}.
\end{align}
The approximation sign in equation \eqref{alpharestriction} is in the sense that both sides agree to within an $O(1)$ factor.
Using these simplifications, we now write the reduced propagator as
\begin{equation}
\Delta a_{\bm{\alpha}}^{(m)}\sim \sum_k \Delta^2_k(p_{k\bm{\alpha}}a_{\bm{\alpha}}).
\end{equation}

Here, $\Delta_k^2$ is the second centered difference operator over $\alpha_k$. Let $\mathbf{N}$ be a vector with components $N_k$.
We can write the reduced propagator in generating function form by letting

\begin{equation}
Q^{(m)}(\mathbf{x},\mathbf{y})=\sum_{\bm{\alpha}} a_{\bm{\alpha}}^{(m)} \mathbf{x}^{\bm{\alpha}} \mathbf{y}^{\mathbf{N}-\bm{\alpha}}.
\end{equation}
Here, the vector powers of $\mathbf{x},\mathbf{y}$ are in the sense of Laurent Schwartz, i.e., $\mathbf{x}^{\bm{\alpha}}=\prod_i x_i^{\alpha_i}$ \cite{john}. We can write the reduced propagator in terms of this function and obtain the following equation for the spectral problem:
\begin{equation}
(\lambda-1)G\sim\sum_k \frac{(x_k-y_k)^2}{NN_k}\frac{\partial^2G}{\partial x_k\partial y_k}.
\end{equation}
Here, $G$ takes the same form as $Q^{(m)}$, and the coefficients are approximate values of the eigenvectors as we had done before.
To solve this, let $u_k=x_k-y_k$ and $G(\mathbf{x},\mathbf{y})=H(\mathbf{u},\mathbf{y})$ to obtain
\begin{equation}
(\lambda-1)H\sim\sum_k \frac{u_k^2}{NN_k}\left(\frac{\partial^2H}{\partial u_k\partial y_k}-\frac{\partial^2H}{\partial u_k^2}\right).
\end{equation}
This translates to the finite difference equation for the coefficients of $H$:

\begin{multline}
(\lambda-1)b_{\bm\alpha}\sim\sum_k \bigg[\frac{(\alpha_k-1)(N_k-\alpha_k+1)b_{\bm{\alpha}-k}}{NN_k}\\
-\frac{\alpha_k(\alpha_k-1)b_{\bm{\alpha}}}{NN_k}\bigg].\\
\end{multline}

By requiring that there is a singularity in the difference equation,
we find a general form of the spectrum as
\begin{align}
\lambda&\sim 1-\sum_k\frac{\alpha_k(\alpha_k-1)}{NN_k}\\
&\approx 1-\sum_k \left(\frac{\alpha_k}{N_k}\right)^2\frac{N_k}{N}.\label{approxspectrum}
\end{align}

To proceed, we need to apply the restriction on $\alpha_k/N_k$ given in equation \eqref{alpharestriction}. In particular, we are looking for small values of $\alpha_k$, which will maximize $\lambda$ such that $\lambda <1$. Notice that when $\alpha_k$ is small and non-zero, we find from equation \eqref{alpharestriction} that

\begin{equation}
\frac{\alpha_k}{N_k}= O\left(\frac{k}{N\mu_1}\right)
\end{equation}

Use this in equation \eqref{approxspectrum} to find that
\begin{align}
&\lambda\approx 1-O\left(\frac{1}{N^2\mu_1^2}\sum_k k^2 \frac{N_k}{N}\right)\\\
&1-\lambda=O\left(\frac{\mu_2}{N^2\mu_1^2}\right).
\end{align}\\

This gives an asymptotic approximation of the spectral gap. The calculation for
finding the moments of consensus time given in section \ref{sec:applications} can be applied here to attain the result given in equation \eqref{moments}. For the expected time to consensus, this result states that $E[T]=O\left(\frac{\mu_1^2}{\mu_2}N^2\right)$, which is consistent with known results regarding consensus times on uncorrelated heterogeneous networks \cite{sood, vazquez}.

\begin{figure}[h!]
\includegraphics[scale=.5]{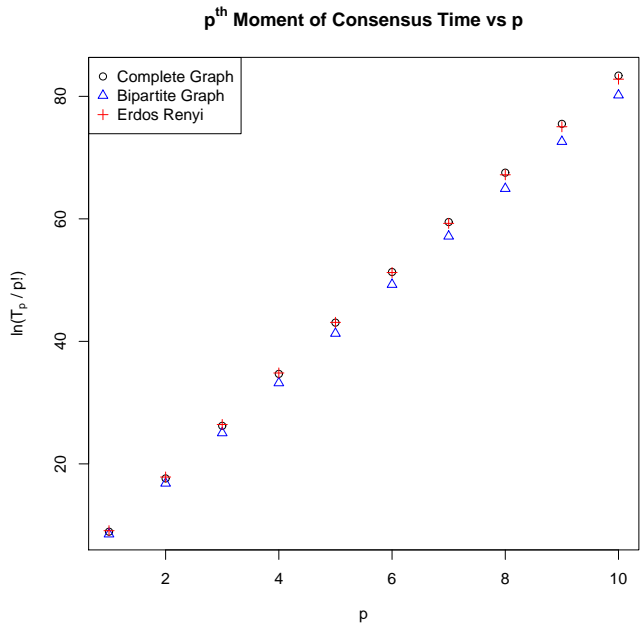}
\caption{Monte Carlo simulation of the voter model on the complete graph, complete bipartite graph, and uncorrelated networks. We fix $N=100$ and let $p$ range from $1...10$.}\label{montecarlo1}
\end{figure}

\begin{figure}[h!]
\includegraphics[scale=.5]{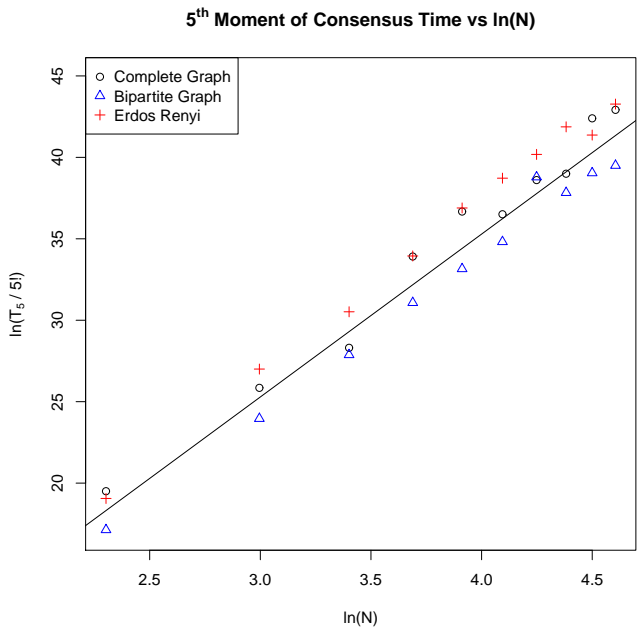}
\caption{Monte Carlo simulation of the voter model on the complete graph, complete bipartite graph, and Erd\H{o}s R\'{e}nyi networks. In this simulation, we fix $p=5$ and let $N$ range from $10,20,...,100$. In addition, a line of slope $2p=10$ is plotted, which is the prediction given by the theory.}\label{montecarlo2}
\end{figure}

For each network that we consider, we test the analytical result for the moments of consensus times against
Monte Carlo simulation data. To this end, we conduct two simulations. In the first simulation, we take $N=100$, and run
the model $100$ times for each network we have considered. For the bipartite graph, we fix $N_1=4N_2$ in all simulations.
For the uncorrelated network, we use the Erd\H{o}s R\'{e}nyi model with the probability of linking two nodes set to be $5/N$ \cite{erdos}.
We let $T_p$ be the numerical result for the $pth$ moment of consensus time for the complete and bipartite simulations.
For the Erd\H{o}s R\'{e}nyi networks, we normalize the time to consensus by $\mu_1^2/\mu_2$ prior to computing $T_p$.
We do this since equation \eqref{moments} predicts that there is a linear relationship between $\ln(T_p/p!)$ and $p$ in each network.
In figure \ref{montecarlo1}, we see that the simulation provides a linear agreement between these quantities.

The second simulation that we provide tests the asymptotic relationship between an arbitrary moment of consensus time with $N$.
For each network, we simulate the voter model for $N=10,20,...,100$. For each of these values of $N$, we run the model $100$ times
and numerically compute the $5th$ moment in each case. This is to test the dependence of the moments of consensus time with $N$.
In equation \eqref{moments}, we expect to see that there is a linear relationship between $\ln(T_p/p!)$ and $\ln N$.
Furthermore, since we take $p=5$ for each $N$, we expect that the slope should be about $10$ in each case. Figure \ref{montecarlo2}
shows that the trend is in fact linear and that the prediction about the slope is accurate.

\section{Conclusions}
\label{sec:conclusions} We have successfully derived analytical
solutions to the voter model for several different networks. In
particular, the corresponding spectral problem can be solved exactly
for the complete graph, and approximately for the other networks.
This allows us to estimate the time $m$ probability distribution of
the model. In addition to the expected time to consensus, we were
able to extract all moments thereof by utilizing the spectral
analysis we have presented here. The expected time to consensus can
also be found by using first step analysis \cite{zhang,sood}, but
that method cannot find any higher moments. Using spectral analysis
to this end is a significant improvement on the theory.

For the complete graph, the complete set of eigenvalues and
eigenvectors had been found. With this, we found the solution to the
corresponding differential eigenvalue problem. The solutions of
which are found to be hypergeometric functions that have terminating
series expressions. Exact formulae for the expected time to
consensus and the expected frequency of each macrostate prior to
consensus are also given, though there are other quantities of
interest that can be found using the theory given above. This can be
done since we can make use of the detailed information given to us
by the eigenvectors. For the other networks, the region of the plane
we have considered had been reduced in order to solve it. As such,
we expect that the results we obtain for the eigenvectors will be
approximate.

Since the procedures are easily generalizable,
there is great potential for applying these techniques to other problems. This includes extensions of the voter model
and possibly other models, such as invasion processes or link dynamics models. Studying these models with the spectral analysis presented here could provide detailed solutions such as the moments of consensus time, local times, and the $m$ step propagator for these other models as well.

\section*{Acknowledgement}

This work was supported in part by the Army Research Office Grant
No. W911NF-09-1-0254 and W911NF-12-1- 467 0546. The views and
conclusions contained in this document are those of the authors and
should not be interpreted as representing the official policies,
either expressed or implied, of the Army Research Office or the U.S.
Government.


\bibliography{Voter14b.bib}{}

\end{document}